\input ./style/arxiv-general.cfg
\documentclass[MSNbibl,number,citesort,seceqn,dvips]{arxbj}
\makeatletter
   \@ifpackageloaded{graphicx}{}{\usepackage{graphicx}}
\makeatother
\usepackage{mathbh}

\volume{22}
\issue{3}
\pubyear{2016}
\firstpage{1520}
\lastpage{1534}
\doi{10.3150/15-BEJ701}
\docsubty{FLA}

\makeatletter
\newremark{nota}{Notation}
\newtheorem{TTheorem}{Theorem}
\newtheorem{Theorem}{Theorem}[section]
\newtheorem{Lemma}[Theorem]{Lemma}
\newtheorem{Proposition}[Theorem]{Proposition}
\newtheorem{Corollary}[Theorem]{Corollary}
\newremark{Remark}[Theorem]{Remark}
\makeatother

\begin{document}
\begin{frontmatter}

\title{Performance of empirical risk minimization in linear aggregation}
\runtitle{ERM in linear aggregation}

\begin{aug}
\author[A]{\inits{G.}\fnms{Guillaume}~\snm{Lecu\'e}\corref
{}\thanksref
{A}\ead
[label=e1]{guillaume.lecue@cmap.polytechnique.fr}}
\and
\author[B]{\inits{S.}\fnms{Shahar}~\snm{Mendelson}\thanksref{B}\ead
[label=e2]{shahar@tx.technion.ac.il}}
\address[A]{CNRS, CMAP, Ecole Polytechnique, 91120 Palaiseau, France.\\
\printead{e1}}
\address[B]{Department of Mathematics, Technion, I.I.T, Haifa 32000,
Israel. \printead{e2}}
\end{aug}

%
\received{\smonth{3} \syear{2014}}
%
\revised{\smonth{2} \syear{2015}}

%
\begin{abstract}
We study conditions under which, given a dictionary $F=\{f_1,\ldots
,f_M\}$
and an i.i.d. sample $(X_i,Y_i)_{i=1}^N$,
the empirical minimizer in $\operatorname{span}(F)$ relative to the squared
loss, satisfies that with high probability
\[
R \bigl(\tilde f^{\mathrm{ERM}} \bigr)\leq\inf_{f\in\operatorname
{span}(F)}R(f)+r_N(M),
\]
where $R(\cdot)$ is the squared risk and $r_N(M)$ is of the order of $M/N$.

Among other results, we prove that a uniform small-ball
estimate for functions in $\operatorname{span}(F)$ is enough to
achieve that
goal when the noise is
independent of the design.
\end{abstract}

%
\end{frontmatter}

\section{Introduction and main results} \label{secintro}
Let $(\mathcal{X},\mu)$ be a probability space, set $X$ to be distributed
according to $\mu$ and put $Y$ to be an unknown target random variable.

In the usual setup in learning theory, one observes $N$ independent
couples $(X_i,Y_i)_{i=1}^N$ in $\mathcal{X}\times\mathbb{R}$,
distributed according
to the joint distribution of $X$ and $Y$. The goal is to construct a
real-valued function $f$ which is a good guess/prediction of $Y$. A
standard way of measuring the prediction capability of $f$ is via the
risk $R(f)=\mathbb{E}(Y-f(X) )^2$. The conditional expectation
\[
R(\hat f)=\mathbb{E} \bigl( \bigl(Y-\hat f(X) \bigr)^2|(X_i,Y_i)_{i=1}^N
\bigr)
\]
is the risk of the function $\hat f$ that is chosen by the procedure,
using the observations $(X_i,Y_i)_{i=1}^N$.

There are many different ways in which one may construct learning
procedures (see, e.g., the books \cite{DGL96,MR1741038,MR2450103,MR726392,MR2722294,MR2724359} for numerous examples), but in general, there is no
`universal' choice of an optimal learning procedure.

The variety of learning algorithms motivated the introduction of
aggregation or ensemble methods, in which one combines a batch or
\textit{dictionary}, created by learning procedures, in the hope of
obtaining a
function with `better' prediction capabilities than individual members
of the dictionary.

Aggregation procedures have been studied extensively (see, e.g., \cite
{MR1792783,MR1775638,MR2163920,MR2044592,MR1946426,MR1762904,MR2458184,DT07,TsyCOLT07}
and references therein), and among the more well-known aggregation
procedures are boosting \cite{MR2920188} and bagging \cite{MR726392}.

Our aim is to explore the problem of \textit{linear aggregation}:
given a
dictionary $F=\{f_1,\ldots,f_M\}$, one wishes to construct a procedure
$\tilde f$ whose risk is almost as small as the risk of the best
element in the linear span of the dictionary, denoted by $\operatorname
{span}(F)$; namely, a procedure which ensures that with high probability
\begin{equation}
\label{eqoracleineq} R(\tilde f)\leq\inf_{f\in\operatorname
{span}(F)}R(f)+r_N(M).
\end{equation}
This type of inequality is called an \textit{oracle inequality} and the
function $f^*$ for which $R(f^*)=\inf_{f\in\operatorname
{span}(F)}R(f)$ is
called the \textit{oracle}.

Of course, in (\ref{eqoracleineq}) one is looking for the smallest
possible residual term $r_N(M)$, that holds uniformly for all choices
of couples $(X,Y)$ and dictionaries $F$ that satisfy certain assumptions.

The linear aggregation problem has been studied in \cite{MR1775638} in
the Gaussian white noise model; in \cite{TsyCOLT07,MR2351101} for the
Gaussian model with random design; in \cite{MR2356821} for the density
estimation problem and in \cite{MR2906886} in the learning theory
setup, under moment conditions. And, based on these cases, it appears
that the best possible residual term $r_N(M)$ that one may hope for is
of the order of $M/N$.

This rate is usually called the \textit{optimal rate of linear
aggregation} and, in fact, its optimality holds in some minimax sense,
introduced in \cite{TsyCOLT07}.

The only procedure we will focus on here is empirical risk
minimization (ERM) performed in the span of the dictionary:
\[
\hat f^{\mathrm{ERM}}\in\mathop{\operatorname{arg}\operatorname{min}}_{f\in
\operatorname{span}(F)}
R_N(f)\qquad\mbox{where } R_N(f)= \frac{1}{N}\sum
_{i=1}^N \bigl(Y_i-f(X_i)
\bigr)^2.
\]
We do not claim that ERM
is always the best procedure for the linear aggregation problem, but
rather, our aim is to identify
conditions under which it achieves the optimal rate of $M/N$.

The benchmark result on the performance of ERM in linear aggregation is
Theorem~2.2 in \cite{MR2906886}. To
formulate it, let $F$ be a
dictionary of cardinality $M$ and set $f^*$ to be the oracle in
$\operatorname{span}(F)$ (i.e., $R(f^*)=\inf_{f\in\operatorname
{span}(F)}R(f)$). We also
denote by $L_p$ for $1\leq p\leq\infty$ the Banach spaces
$L_p(\mathcal{X},
\mu)$, and in particular, $\Vert f\Vert_{L_2}= (\mathbb{E}f(X)^2 )^{1/2}$.

\begin{Theorem}[\cite{MR2906886}]\label{theocatoni}
Assume that $\mathbb{E}(Y-f^*(X))^4< \infty$ and that for every $f\in
\operatorname{span}(F)$,
\begin{equation}
\label{eqcatoni} \Vert f\Vert_{L_\infty}\leq\sqrt{B} \Vert f
\Vert_{L_2}.
\end{equation}
If $x>0$ satisfies that $2/N\leq2\exp(-x)\leq1$ and
\[
N\geq1280B^2 \biggl[3BM + x+\frac{16B^2 M^2}{N} \biggr],
\]
then with probability at least $1-2\exp(-x)$,
\[
R \bigl(\hat f^{\mathrm{ERM}} \bigr)-R \bigl(f^* \bigr)\leq1920 B\sqrt{\mathbb{E}
\bigl(Y-f^*(X) \bigr)^4} \biggl[\frac{3BM+x}{N}+\frac{16 B^2
M^2}{N^2}
\biggr].
\]
\end{Theorem}

It follows from Theorem~\ref{theocatoni} that under an $L_4$
assumption on $Y-f^*(X)$ and the equivalence between the $L_2$ and
$L_\infty$ norms on the span of $F$, ERM achieves a rate of
convergence of order $B^2M/N$ when $N \geq c B^3M$ for an absolute
constant $c$.

However, it should be noted that the best probability estimate one may
obtain in Theorem~\ref{theocatoni}
is $1-2/N$; also, it is possible to show that the constant $B$ defined
in (\ref{eqcatoni}) is \textit{necessarily larger than the dimension $M$
of $\operatorname{span}(F)$}. For the sake of completeness, we shall
provide a
proof of that fact in the \hyperref[secappendix]{Appendix}. Therefore, the rate that
Theorem~\ref{theocatoni} guarantees is, at best, of the order of
$M^3/N$, to achieve that rate, at least $N \geq cM^4$ observations are
needed, and even with that sample size, the probability estimate is, at
best, $1-2/N$. This estimate is far from the anticipated rate of $M/N$,
which should be achieved
when $N \geq c M$ and preferably, with significantly higher probability.

Nevertheless, the optimal rate of $M/N$ can be obtained by relaxing
assumption~(\ref{eqcatoni}) and using a different method of proof.
Recall that the $\psi_2$ norm of a
function $f$ is
\[
\Vert f\Vert_{\psi_2}=\inf \bigl\{C>0 \dvt\mathbb{E}\exp
\bigl(f^2(X)/C^2 \bigr)\leq2 \bigr\}.
\]
One may show that $\Vert f\Vert_{\psi_2}\leq c
\Vert f\Vert_{L_\infty}$ for a suitable absolute constant $c$ (see, e.g.,
Section~1 in \cite{MR3113826}). Therefore, assuming that the $\psi
_2$-norm and the $L_2$-norm are equivalent in $\operatorname{span}(F)$
is a weaker
requirement than the one in (\ref{eqcatoni}). The assumption that for
every $f\in\operatorname{span}(F)$,
\begin{equation}
\label{eqpsi-2-L-2} \Vert f\Vert_{\psi_2}\leq\sqrt{C}\Vert f
\Vert_{L_2},
\end{equation}
means that $\operatorname{span}(F)$ is a \textit{sub-Gaussian class},
following the definition from \cite{LM13}. To put this assumption in
some perspective, there are numerous examples of sub-Gaussian classes
(the simplest of which are classes of linear functionals on $\mathbb{R}^M$
endowed with a sub-Gaussian design) for which the equivalence constant
$C$ is an absolute constant, unlike the constant $B$ in (\ref
{eqcatoni}), which is at least $M$.

Naturally, the analysis of ERM under a sub-Gaussian assumption requires
a more sophisticated technical machinery than in situations in which
the $L_2/L_\infty$ equivalence assumption used in Theorem~\ref
{theocatoni} holds. Invoking the main result from \cite{LM13}, one
can show that if $Y-f^*(X)$ is
sub-Gaussian and $\operatorname{span}(F)$ is a \textit{sub-Gaussian class},
then for every $x>0$,
ERM achieves a rate $r_N(M)=c_1xM/N$ with probability at least $1-\exp
(-c_2xM)$.

Although the sub-Gaussian case is interesting, the goal of this note is
the study of ERM as a linear aggregation procedure under much weaker
assumptions.

\renewcommand{\theTTheorem}{A}
\begin{TTheorem}\label{A}
Let $F=\{f_1,\ldots,f_M\}$ and assume that
there are constants $\kappa_0$ and $\beta_0$ for which
\begin{equation}
\label{eqsmallball} P \bigl\{ \bigl\vert f(X) \bigr\vert\geq\kappa_0 \Vert f
\Vert_{L_2} \bigr\}\geq\beta_0
\end{equation}
for every $f\in\operatorname{span}(F)$. Let $N \geq(400)^2M/\beta
_0^2$ and
set $\zeta=Y-f^*(X)$. Assume further that one of the following two
conditions holds:
\begin{enumerate}[2.]
\item[1.] $\zeta$ is independent of $X$ and $\mathbb{E}\zeta^2\leq
\sigma^2$, or
\item[2.] $|\zeta|\leq\sigma$ almost surely.
\end{enumerate}
Then, for every $x>0$, with probability at least $1-\exp(-\beta_0^2
N/4)-(1/x)$,
\[
\bigl\Vert\hat{f}^{\mathrm{ERM}}-f^*\bigr\Vert_{L_2}^2=R \bigl(\hat
f^{\mathrm{ERM}} \bigr)-\min_{f\in\operatorname{span}(F)}R(f)\leq \biggl(
\frac{16}{\beta_0 \kappa_0^2} \biggr)^2\frac
{\sigma^2 Mx}{N}.
\]
\end{TTheorem}

Since the loss is the squared one, one has to assume that $Y$ and
functions in $\operatorname{span}(F)$ have a second moment. It follows from
Theorem~\ref{A} that in some cases, this is (almost) all that is needed for
an optimal rate. Indeed, if $\zeta=Y-f^*(X)$ is independent
of the design $X$ -- as is the case in any regression model with
independent noise $Y=f^*(X)+\zeta$, and if (\ref{eqsmallball})
holds, ERM achieves the optimal rate $M/N$.

\begin{Corollary}\label{cororegression-model}
Consider the regression model $Y=f^*(X)+\zeta$ where $\zeta$ is a
mean-zero noise that is independent of $X$. Assume that $\zeta\in L_2$ and
that $f^*\in\operatorname{span}(F)$. If $\operatorname{span}(F)$
satisfies (\ref{eqsmallball}) and $N \geq(400)^2M/\beta_0^2$, then
for every $x>0$, with
probability at least $1-\exp(-\beta_0^2N/4)-1/x$,
\[
\bigl\Vert\hat{f}^{\mathrm{ERM}}-f^*\bigr\Vert_{L_2}^2\leq \biggl(
\frac{16}{\beta_0 \kappa
_0^2} \biggr)^2\frac{\sigma^2 Mx}{N}.
\]
\end{Corollary}

From a statistical point of view, (\ref{eqsmallball}), which is a
\textit{small-ball assumption} on $\operatorname{span}(F)$, is a quantified
version of \textit{identifiability}. Indeed,
consider the statistical model $\mathcal{M}= \{\mathbb{P}_f \dvt f\in
\operatorname{span}(F)\}$
where $\mathbb{P}_f$ is the probability distribution of the couple $(X,Y)$,
$Y=f(X)+\zeta$ and $\zeta$ is, for instance, a Gaussian noise that is
independent of $X$. Assuming that $\mathcal{M}$ is identifiable is
equivalent to
having $P (|f(X)-g(X)|>0 )>0$ for every $f,g\in\operatorname{span}(F)$,
which, by linearity, is equivalent to
$P (|f(X)|>0 )>0$ for every $f\in\operatorname{span}(F)$. Comparing
this with the small-ball
condition in~(\ref{eqsmallball}) shows that the latter is just a
`robust' version of identifiability.

It is possible to slightly modify the assumptions of Theorem~\ref{A} and
still obtain the same type of estimate. For example, it is
straightforward to verify that the small-ball condition (\ref
{eqsmallball}) holds when the $L_2$ and $L_p$ norms are equivalent on
$\operatorname{span}(F)$ for some $p>2$. This type of
$L_p/L_2$ equivalence assumption on $\operatorname{span}(F)$ is weaker
than the
equivalence between the $L_{\psi_2}$ and the $L_2$ norms in (\ref
{eqpsi-2-L-2}) because for every $p\geq1$,
$\Vert f\Vert_{L_p} \leq c \sqrt{p}\Vert f\Vert_{\psi_2}$ for a suitable
absolute constant $c$. And, it is clearly weaker than the $L_\infty
/L_2$ equivalence assumption~(\ref{eqcatoni}) used in Theorem~\ref
{theocatoni}.

It turns out that if the $L_2$ and $L_4$ norms are equivalent on
$\operatorname{span}(F)$, one may obtain the optimal rate for an
arbitrary target $Y$,
as long as
$\zeta=Y-f^*(X)$ has a fourth moment. The difference between such a
result and
Theorem~\ref{A} is that $\zeta$ need not be independent of $X$, nor must it
be bounded.

\begin{Theorem}\label{coroL4-L2}
There exist absolute constants $c_0,c_1$ and $c_2$ for which the
following holds. Assume that there exists $\theta_0$ for which
\begin{equation}
\label{eqweakmoment} \Vert f\Vert_{L_4}\leq\theta_0 \Vert f
\Vert_{L_2}
\end{equation}
for every $f\in\operatorname{span}(F)$, and let $N \geq(c_0\theta
_0^4 )^2 M$. Set $\zeta=Y-f^*(X)$ and put $\sigma= (\mathbb{E}\zeta
^4 )^{1/4}$. Then, for every $x>0$, with probability at least
$1-\exp(-N/(c_1\theta_0^8))-(1/x)$,
\[
\bigl\Vert\hat f-f^*\bigr\Vert_{L_2}^2=R(\hat f)-\min
_{f\in\operatorname{span}(F)}R(f)\leq c_2\theta_0^{12}
\cdot\frac{\sigma^2 Mx}{N}.
\]
\end{Theorem}

\begin{Remark}
One may show that a possible choice of constants in Theorem~\ref
{coroL4-L2} is $c_0=1600$, $c_1=64$ and $c_2=(256)^2$, but since we
have not made any real attempt of optimizing the choice of constants --
because identifying the correct rate is the main focus of this note --
we will not keep track of the values of constants in what follows.
\end{Remark}

One example in which Theorem~\ref{coroL4-L2} may be used is the
regression problem with a misspecified model: $Y=f_0(X)+W$ where
the regression function $f_0$ may not be in the model $\operatorname{span}(F)$
and $\zeta=(f_0-f^*)(X)+W$ has a fourth moment. If $\operatorname{span}(F)$
satisfies~(\ref{eqsmallball}), then with high probability,
\begin{equation}
\label{eqmisspecification} \bigl\Vert\hat{f}-f^* \bigr\Vert_{L_2}^2= \bigl\Vert\hat
f-f_0\bigl\Vert_{L_2}^2-\bigl\Vert f_0-f^*
\bigr\Vert_{L_2}^2 \leq c(\theta_0) \bigl(\mathbb{E}
\zeta^4 \bigr)^{1/2}\frac{M}{N},
\end{equation}
for a constant $c(\theta_0)$ that only depends on $\theta_0$. Hence,
one may select $M$ as the
solution of an optimal trade-off between the variance term $ (\mathbb{E}
\zeta^4 )^{1/2}M/N$ and
the bias; we refer the reader to Chapter~1 in \cite{MR2724359} for
techniques of a similar flavour.

The standard way of analyzing the performance of ERM is via certain
trade-offs between concentration and complexity. However, in the case
we study here, the functions involved may have `heavy tails', and
empirical means do not exhibit strong, two-sided concentration around
their true means -- which is a crucial component in the standard method
of analysis. Therefore, a completely different path must be taken if
one is to obtain the results formulated above.

The method we shall employ here has been introduced in \cite
{Shahar-COLT,Shahar-LWCG} for problems in Learning Theory;
in \cite{Shahar-Gelfand} in the context of the geometry of convex
bodies; in \cite{Shahar-Vladimir} for applications in
random matrix theory; and in \cite{LMcompressed} for Compressed Sensing.

Obviously, and regardless of the method of analysis, the (seemingly)
unsatisfactory probability estimate is the price one pays for the
moment assumptions on the `noise' $Y-f^*(X)$. The next result shows
that without stronger moment assumptions, only weak polynomial
probability estimates are true.

\begin{Proposition}\label{propexpo-moments}
Let $x\geq1$, assume that $N\geq c_0M$ for a suitable absolute
constant $c_0$ and that $X$ is the standard Gaussian vector in $\mathbb{R}^M$.
There exists a mean-zero, variance one random variable~$\zeta$, that
is independent of $X$ and for which the following holds.

Fix $t^* \in\mathbb{R}^M$ and consider the model $Y= \langle
X,t^*  \rangle+\zeta$.
With probability at least $c_1/x$, ERM produces $\hat t\in
\operatorname{arg}\operatorname{min}
_{t\in\mathbb{R}^M}\sum_{i=1}^N (Y_i- \langle X_i,t
\rangle )^2$
that satisfies
\[
\bigl\Vert\hat t-t^*\bigr\Vert_2^2=R(\hat t)-R \bigl(t^* \bigr)
\geq \frac{c_2xM}{N},
\]
where $c_1$ and $c_2$ are absolute constants and
$R(t)=\mathbb{E}(Y- \langle X,t  \rangle )^2$ is the
squared risk of $t$.
\end{Proposition}

Note that the class of linear functional $\{ \langle\cdot,t
 \rangle \dvt t\in\mathbb{R}^M\}
$ is a linear space of dimension $M$
and it satisfies the small-ball condition when $X$ is the standard
Gaussian vector (actually, this
class is sub-Gaussian). It follows from Proposition~\ref
{propexpo-moments} that there is no hope
of obtaining an exponential probability bound on the excess risk of ERM
under an $L_2$-moment assumption on the noise -- only polynomial bounds
are possible. In particular, the probability estimate
obtained in Theorem~\ref{A} under the $L_2$-assumption on the noise cannot be
improved.

Finally, we would like to address the problem of linear aggregation
under the classical boundedness assumptions: that $|Y|\leq1$ and
$|f(X)|\leq1$ almost surely for every $f \in F$.

These are the standard assumptions that have been considered for
the three problems of aggregation with a random
design. For instance, optimal rates of aggregation have been obtained under
these assumptions for the model selection aggregation problem in
\cite{MR2529440,Audibert1,LecRig} and for the convex aggregation problem
in \cite{MR3160549}. And, it has been established that while ERM is
suboptimal for the model
selection aggregation problem (see, e.g., Section~3.5 in \cite
{MR2163920} or
\cite{LM2}), it is optimal for the convex aggregation problem.
However, the optimality of ERM in the linear aggregation problem under
the boundedness
assumption was left open. The final result of this article addresses
that problem -- and it turns out that the answer is negative in a
very strong way.

\begin{Proposition}\label{propsub-optimality-ERM-linear-agg}
For every $0<\eta<1$ and integers $N$ and $M$, there exists a couple
$(X,Y)$ and a dictionary $F=\{f_1,\ldots,f_M\}$ with the following
properties:
\begin{enumerate}[2.]
\item[1.] $|Y|\leq1$ almost surely and $|f(X)|\leq1$ almost surely
for every $f\in F$.
\item[2.] With probability at least $\eta$, for every $\kappa>0$
there is some
\[
\hat f^{\mathrm{ERM}}\in\mathop{\operatorname{arg}\operatorname{min}}_{f\in
\operatorname{span}(F)}
\frac{1}{N}\sum_{i=1}^N
\bigl(Y_i-f(X_i) \bigr)^2
\]
for which\vspace*{-4pt}
\[
R \bigl(\hat f^{\mathrm{ERM}} \bigr)\geq\inf_{f\in\operatorname
{span}(F)}R(f)+\kappa.
\]
\end{enumerate}
\end{Proposition}

Proposition~\ref{propsub-optimality-ERM-linear-agg} shows that even
if one assumes that $|Y|\leq1$ and $|f(X)|\leq1$
almost surely for every function in the dictionary, and despite the
convexity of $\operatorname{span}(F)$, the empirical risk minimization
procedure
performs poorly. This illustrates the major difference between assuming
that the class is well bounded in $L_\infty$ and assuming that the
$L_2$ and $L_p$ norms are equivalent on its span: while the latter
suffices for an optimal bound, the former is rather useless.

An obvious outcome of Proposition~\ref
{propsub-optimality-ERM-linear-agg} is that ERM
should not be used to solve the linear aggregation problem under the
boundedness assumption and one has to look for different procedures in
the bounded setup. It should also be noted that since Proposition~\ref
{propsub-optimality-ERM-linear-agg} is a non-asymptotic lower bound
and $X$ may depend on $N$ and $M$, the
asymptotic result appearing in Theorem~2.1 in \cite{MR2906886} does
not apply here.

\begin{nota*}
For every function $f$, let
$\Vert f\Vert_{L_p}= (\mathbb{E}|f(X)|^p )^{1/p}$. The excess loss of a
function $f\in\operatorname{span}(F)$ is defined for every $x\in
\mathcal{X}$ and
$y\in\mathbb{R}$ by
\[
\mathcal{L}_f(x,y)= \bigl(y-f(x) \bigr)^2-
\bigl(y-f^*(x) \bigr)^2;
\]
thus, $R(f)-R(f^*)=P\mathcal{L}(X,Y)\geq0$. The empirical
measure over the data is denoted by $P_N$ and\vspace*{-4pt}
\[
P_N\mathcal{L}_f=\frac{1}{N}\sum
_{i=1}^N \bigl(Y_i-f(X_i)
\bigr)^2- \bigl(Y_i-f^*(X_i)
\bigr)^2.
\]
For every\vadjust{\goodbreak} vector $x\in\mathbb{R}^M$, let
$\Vert x\Vert_{\ell_p^M}= (\sum_{j=1}^M |x_j|^p )^{1/p}$ be its
$\ell_p^M$-norm.

Finally, all absolute constants are denoted by $c_1,c_2$, etc. Their
value may change from line to line. We write $A \lesssim B$ if there is
an absolute constant $c$ for which $A \leq cB$, and $A \lesssim_\alpha
B$ if $A \leq c(\alpha)B$ for a constant $c$ that depends only on
$\alpha$.
\end{nota*}

\section{Proofs of Theorem~\texorpdfstring{\protect\ref{A}}{A} and Theorem~\texorpdfstring{\protect\ref{coroL4-L2}}{1.3}}
\label{sectheoB}
The starting point of the proof of Theorem~\ref{A} is the same as in \cite
{LM13,Shahar-COLT,LMphaseretreival,Shahar-LWCG}: a decomposition of
the excess loss function
\begin{equation}
\label{eqdecompositionloss} \mathcal{L}_f(x,y)= \bigl(f^*(x)-f(x)
\bigr)^2+2 \bigl(y-f^*(x) \bigr) \bigl(f^*(x)-f(x) \bigr)
\end{equation}
to a sum of quadratic and linear terms in $(f-f^*)(X)$. The idea of the
proof is to control the quadratic term from below using a
`small-ball' argument, and the linear term from above using standard
methods from empirical processes theory. A combination of these two
bounds suffices to show that if $\|f-f^*\|_{L_2} \geq r_N^*$ for an
appropriate choice of $r_N^*$, the quadratic term dominates the
linear one, and in particular, for such functions $P_N \mathcal
{L}_f>0$. Since
the empirical excess loss of the empirical minimizer is non-positive,
it follows that $\|\hat{f}-f^*\|_{L_2} < r_N^*$.

\begin{Lemma}\label{lemsmallball}
There exists an absolute constant $c_0$ for which the following holds.
Assume that there are $\kappa_0$ and $\beta_0$ for which
\[
P \bigl( \bigl\vert f(X) \bigr\vert\geq\kappa_0 \Vert f\Vert_{L_2}
\bigr)\geq\beta_0
\]
for every $f\in\operatorname{span}(F)$. If $N \geq c_0M/\beta_0^2$,
then with
probability at least $1-\exp(-\beta_0^2 N/4)$, for every $f\in
\operatorname{span}(F)$,
\[
\bigl\vert \bigl\{i\in\{1,\ldots,N\}\dvt \bigl\vert f(X_i) \bigr\vert\geq
\kappa_0 \Vert f\Vert_{L_2} \bigr\} \bigr\vert \geq
\frac{\beta_0 N}{2}.
\]
\end{Lemma}

\begin{pf}
Let $x>0$ and set
\[
H= \sup_{f\in\operatorname{span}(F)} \Biggl\vert\frac{1}{N}\sum
_{i=1}^N \mathbh{1}_{ \{\vert f(X_i)\vert\geq\kappa_0 \Vert f\Vert
_{L_2} \}}-P \bigl( \bigl\vert
f(X) \bigr\vert\geq\kappa_0 \Vert f\Vert_{L_2} \bigr) \Biggr\vert.
\]
Set $W=(f_1(X),\ldots,f_M(X))$ -- a random vector endowed on $\mathbb
{R}^M$ by
the dictionary $F$ and the
random variable $X$. Note that $\operatorname{span}(F)= \{\sum_{j=1}^M t_j f_j \dvt(t_1,\ldots,t_M)\in\mathbb{R}^M \}$
and set $\|t\|_{L_2}=\|\sum_{j=1}^M t_j f_j\|_{L_2}$.

Since $N$ independent copies of $X$, $X_1,\ldots,X_N$, endow $N$
independent copies
of $W$, denoted by $W_1,\ldots,W_N$, it follows that
\[
H=\sup_{t \in\mathbb{R}^M } \Biggl\vert\frac{1}{N}\sum
_{i=1}^N \mathbh{1}_{ \{
\vert \langle t,\cdot  \rangle\vert\geq\kappa_0 \Vert
t\Vert_{L_2} \}
}(W_i)-P
\bigl(\bigl\vert \langle t,W \rangle\bigr\vert\geq\kappa_0 \Vert t
\Vert_{L_2} \bigr) \Biggr\vert.
\]

By the bounded differences inequality (see, e.g., Theorem~6.2 in \cite
{BouLugMass13}), with probability at least $1-\exp(-x^2/2)$,
\begin{equation}
\label{eqmainlem} H\leq\mathbb{E}H+\frac{1}{2}\sqrt{\frac{x}{N}},
\end{equation}
and a standard argument based on the VC-dimension of half-spaces in
$\mathbb{R}
^M$ shows that
\[
\mathbb{E}H = \mathbb{E}H(X_1,\ldots,X_N)\leq
c_1\sqrt{ \frac{M}{N}}
\]
(one may show the $c_1 \leq100$ using a rough estimate on Dudley's
entropy integral combined with
Exercise~2.6.4 in \cite{vanderVaartWellner}). Therefore, if $c_1\sqrt
{M/N}\leq\beta_0/4$ and $(1/2)\sqrt{x/N}= \beta_0/4$, then with
probability at least $1- \exp(-\beta_0^2 N/4)$, $H \leq\beta_0/2$.

Finally,\vspace*{-3pt} since
\[
\inf_{f\in\operatorname{span}(F)} P \bigl( \bigl\vert f(X) \bigr\vert\geq \kappa_0
\Vert f\Vert_{L_2} \bigr)\geq\beta_0
\]
it follows that on the event $\{H\leq\beta_0/2\}$,
\begin{equation}
\label{eqlast1} \inf_{f \in\operatorname{span}(F)} \frac{1}{N}\sum
_{i=1}^N \mathbh{1}_{ \{
\vert f(X_i)\vert\geq\kappa_0 \Vert f\Vert_{L_2} \}}(X_i)
\geq \frac
{\beta_0}{2}.
\end{equation}
Therefore, (\ref{eqlast1}) holds with probability at
least $1- \exp(-\beta_0^2 N/4)$.
\end{pf}

\begin{Lemma}\label{lemlinearprocess}
Let $\zeta=Y-f^*(X)$ and assume that one of the following two
conditions hold:
\begin{enumerate}[2.]
\item[1.] $\zeta$ is independent of $X$ and $\mathbb{E}\zeta^2\leq
\sigma^2$, or
\item[2.] $|\zeta|\leq\sigma$ almost surely.
\end{enumerate}
Then, for every $x>0$, with probability larger than $1-(1/x)$,
\[
\Biggl\vert\frac{1}{N}\sum_{i=1}^N
\bigl(Y_i-f^*(X_i) \bigr) \bigl(f^*(X_i)-f(X_i)
\bigr) \Biggr\vert\leq2\sigma\sqrt{\frac{Mx}{N}} \bigl\Vert f^*-f\bigr\Vert_{L_2}
\]
for every $f\in\operatorname{span}(F)$.
\end{Lemma}

\begin{pf}
Recall that $f^*(X)$ is the best $L_2$-approximation of $Y$ in the
linear space $\operatorname{span}(F)$; hence, $\mathbb
{E}(Y-f^*(X))(f^*(X)-f(X)) = 0$
for every $f\in\operatorname{span}(F)$.

Let $\varepsilon_1,\ldots,\varepsilon_N$ be independent Rademacher
variables that
are also independent of the couples $(X_i,Y_i)_{i=1}^N$. A standard
symmetrization argument shows that
\begin{eqnarray*}
&&\mathbb{E}\sup_{f\in\operatorname{span}(F) \setminus\{f^*\}} \Biggl\vert\frac
{1}{N}\sum
_{i=1}^N \bigl(Y_i-f^*(X_i)
\bigr)\frac
{f^*(X_i)-f(X_i)}{\Vert f^*-f\Vert_{L_2}} \Biggr\vert^2
\\
&&\quad\leq4 \mathbb{E}\sup_{f\in\operatorname{span}(F)\setminus\{f^*\}
} \Biggl\vert\frac
{1}{N}\sum
_{i=1}^N \varepsilon_i
\bigl(Y_i-f^*(X_i) \bigr)\frac
{f^*(X_i)-f(X_i)}{\Vert f^*-f\Vert_{L_2}} \Biggr\vert
^2.
\end{eqnarray*}

Let $T=\{t \in\mathbb{R}^M \dvt\|\sum_{j=1}^M t_j f_j\|_{L_2}=1\}$
and observe
that if $\zeta_1,\ldots,\zeta_N$ are independent copies of $\zeta$, then
\begin{eqnarray*}
&&\mathbb{E}\sup_{f\in\operatorname{span}(F) \setminus\{f^*\}} \Biggl\vert\frac
{1}{N}\sum
_{i=1}^N \varepsilon_i
\bigl(Y_i-f^*(X_i) \bigr)\frac
{f^*(X_i)-f(X_i)}{\Vert f^*-f\Vert_{L_2}}
\Biggr\vert^2
\\
&&\quad= \mathbb{E}\sup_{t \in T} \Biggl\vert\frac{1}{N}\sum
_{i=1}^N \varepsilon_i
\zeta_i \Biggl(\sum_{j=1}^M
t_j f_j(X_i) \Biggr) \Biggr\vert^2=(*).
\end{eqnarray*}
Recall that $W=(f_1(X),\ldots,f_M(X))$ and set $\Sigma$ to be the
covariance matrix associated with~$W$. Let $\Sigma^{-1/2}$ be the
pseudo-inverse of the squared-root of $\Sigma$, set $Z=\Sigma^{-1/2}
W$ and note that $\mathbb{E}\Vert Z\Vert_{\ell_2^M}^2\leq M$.

If $Z_1,\ldots,Z_N$ are independent copies of $Z$, it follows that
\begin{eqnarray*}
(*) &=& \mathbb{E}\sup_{\Vert t\Vert_{\ell_2^M}=1} \Biggl\vert \Biggl\langle t,
\frac {1}{N}\sum_{i=1}^N
\varepsilon_i \zeta_i Z_i \Biggr\rangle
\Biggr\vert^2=\mathbb{E}\Biggl\Vert\frac {1}{N}\sum
_{i=1}^N \varepsilon_i
\zeta_i Z_i\Biggr\Vert_{\ell_2^M}^2
\\
& =& \mathbb{E}\mathbb{E}_{\varepsilon_1,\ldots,\varepsilon_N} \Biggl\Vert\frac{1}{N}\sum
_{i=1}^N \varepsilon_i
\zeta_i Z_i\Biggr\Vert _{\ell_2^M}^2 =
\mathbb{E} \Biggl(\frac{1}{N^2}\sum_{i=1}^N
\zeta_i^2 \Vert Z_i\Vert_{\ell
_2^M}^2
\Biggr)=\frac{\mathbb{E}\zeta^2 \Vert Z\Vert_{\ell_2^M}^2}{N}
\\
& \leq & \frac{\sigma^2\mathbb{E}\Vert Z\Vert_{\ell_2^M}^2}{N},
\end{eqnarray*}
implying that
\begin{eqnarray*}
&&\mathbb{E}\sup_{f\in\operatorname{span}(F)\setminus\{f^*\}} \Biggl\vert\frac{1}{N}\sum
_{i=1}^N \bigl(Y_i-f^*(X_i)
\bigr)\frac{f^*(X_i)-f(X_i)}{\Vert f^*-f\Vert_{L_2}} \Biggr\vert ^2\leq\frac{4\sigma^2 M}{N}.
\end{eqnarray*}
The claim now follows from Markov's inequality.
\end{pf}

\begin{pf*}{Proof of Theorem~\protect\ref{A}}
Combining
Lemma~\ref{lemsmallball} and Lemma~\ref{lemlinearprocess} when
$N \geq c_0M/\beta_0^2$, it follows that with probability at least
$1-\exp(-\beta_0^2N/4)-(1/x)$, if $f\in\operatorname{span}(F)$ and
\begin{equation}
\label{eqlarge-L2} \bigl\Vert\hat f-f^*\bigr\Vert_{L_2}> \frac{16 \sigma
}{\beta_0 \kappa_0^2}\sqrt{
\frac{Mx}{N}},
\end{equation}
one has
\begin{eqnarray*}
&&\frac{1}{N}\sum_{i=1}^N
\bigl(f^*(X_i)-f(X_i) \bigr)^2
\\
&&\quad\geq\kappa_0^2 \bigl\Vert f-f^*\bigr\Vert_{L_2}^2
\bigl| \bigl\{i\dvt \bigl\vert f^*(X_i)-f(X_i) \bigr\vert \geq
\kappa_0\bigl\Vert f-f^*\bigr\Vert_{L_2} \bigr\}\bigr|/N
\\
&&\quad\geq\frac{\beta_0 \kappa_0^2}{2}\bigl\Vert f-f^*\bigr\Vert_{L_2}^2 > 8
\sigma \sqrt{\frac{Mx}{N}}\bigl\Vert f^*-f\bigr\Vert_{L_2}
\\
&&\quad>\frac{2}{N}\sum_{i=1}^N
\bigl(Y_i-f^*(X_i) \bigr) \bigl(f^*(X_i)-f(X_i)
\bigr).
\end{eqnarray*}
Hence, on the same event, if $f\in\operatorname{span}(F)$ and
(\ref{eqlarge-L2}) is satisfied then $P_N\mathcal{L}_f>0$. Since
$P_N\mathcal{L}
_{\hat f^{\mathrm{ERM}}}\leq0$, it follows that
\[
\bigl\Vert\hat f^{\mathrm{ERM}}-f^*\bigr\Vert_{L_2}^2\leq \biggl(
\frac{16 \sigma}{\beta_0
\kappa_0^2} \biggr)^2\frac{Mx}{N}.
\]
\upqed\end{pf*}

\begin{pf*}{Proof of Theorem~\protect\ref{coroL4-L2}}
The proof of
Theorem~\ref{coroL4-L2} is almost identical to the proof of
Theorem~\ref{A}, and we will only outline the minor differences.

The small-ball condition (\ref{eqsmallball}) follows from the
Paley--Zygmund inequality (see, for instance, Proposition~3.3.1 in
\cite{MR1666908}): if $V$ is a real-valued random variable then
\[
P \bigl(\vert V\vert\geq\kappa_0 \bigl(\mathbb{E}V^2
\bigr)^{1/2} \bigr)\geq(1-\kappa_0)^2
\frac{ (\mathbb{E}V^2 )^2}{\mathbb{E}\vert V\vert^4}.
\]
In particular, if $ (\mathbb{E}|V|^4 )^{1/4}\leq\theta_0
(\mathbb{E}|V|^2 )^{1/2}$ then
\[
P \bigl(\vert V\vert\geq(1/2) \bigl(\mathbb{E}V^2
\bigr)^{1/2} \bigr)\geq \bigl(4\theta_0^4
\bigr)^{-1}
\]
and thus the assertion of Lemma~\ref{lemsmallball} holds for $\kappa
_0=1/2$ and $\beta_0=(4\theta_0^4)^{-1}$.

As for the analogous version of Lemma~\ref{lemlinearprocess}, the
one change in its proof is that
\begin{eqnarray*}
&&\mathbb{E}\zeta^2 \Vert Z\Vert_{\ell_2^M}^2\leq
\bigl(\mathbb{E} \zeta^4 \bigr)^{1/2} \bigl(\mathbb{E}\Vert Z
\Vert_{\ell_2^M}^4 \bigr)^{1/2}
\end{eqnarray*}
and
\begin{eqnarray*}
&&\mathbb{E}\Vert Z\Vert_{\ell_2^M}^4= \mathbb{E} \Biggl(\sum
_{j=1}^M \langle e_j,Z
\rangle^2 \Biggr)^2=\mathbb{E}\sum
_{p,q=1}^M \langle e_p,Z
\rangle^2 \langle e_q,Z \rangle^2
\\
&&\quad\leq\sum_{p,q=1}^M \bigl(\mathbb{E}
\langle e_p,Z \rangle^4\mathbb{E} \langle
e_q,Z \rangle^4 \bigr)^{1/2}\leq
\theta_0^4 \sum_{p,q=1}^M
\mathbb{E} \langle e_p,Z \rangle^2\mathbb{E} \langle
e_q,Z \rangle^2=\theta_0^4
M^2.
\end{eqnarray*}
\upqed\end{pf*}

\section{Proof of Proposition~\texorpdfstring{\protect\ref
{propsub-optimality-ERM-linear-agg}}{1.6}} \label{secproof-thm-A}
Fix\vspace*{1pt} $Y=1$ as the target and let $\mathcal{X}=\bigcup_{j=0}^M\mathcal
{X}_j$ be some
partition of $\mathcal{X}$. Consider a random variable $X$ which is
distributed as follows: fix $k\geq M$ to be chosen later; for $1 \leq j
\leq M$, set $P(X\in\mathcal{X}_j)=\frac{1}{k}$ and put $P(X\in
\mathcal{X}
_0)=1-\frac{M}{k}$.

Finally, set
\[
f_j(x)=\cases{ 1, & \quad\mbox{if } $x\in \mathcal{X}_j$,
\vspace*{3pt}
\cr
0, & \quad\mbox{otherwise}}
\]
and put $F=\{f_1,\ldots,f_M\}$.

Note\vspace*{1pt} that $|Y|\leq1$ almost surely and that for every $f\in F$,
$|f(X)|\leq1$ almost surely. It is straightforward to verify that the
oracle in $\operatorname{span}(F)$ is $f^*=\sum_{j=1}^M f_j(\cdot)$,
and thus
\[
\inf_{f\in\operatorname{span}(F)}R(f)=R \bigl(f^* \bigr)=\mathbb {E} \bigl(Y-f^*(X)
\bigr)^2=P(X\in\mathcal{X}_0)=1-\frac{M}{k}.
\]

Let $X_1,\ldots,X_N$ be independent copies of $X$. Given $0< \eta<1$
and $k$ large enough (for instance, $k\geq c(\eta) N/\log M$ for a
sufficiently large constant $c(\eta)$ would suffice), there exists an
event $\Omega_0$ of probability at least
$\eta$ on which the following holds: there exists
$j_0\in\{1,\ldots,M\}$ for which $X_i\notin\mathcal{X}_{j_0}$ for
every $1
\leq i \leq N$ (this is a slight modification of the coupon-collector problem).

For every $j=1,\ldots,M$, let $N_j=|\{i\in\{1,\ldots,N\}\dvt X_i\in
\mathcal{X}
_j\}|$. Hence, for $t\in\mathbb{R}^M$, the empirical risk of $\sum_{i=1}^M
t_j f_j$ is
\begin{eqnarray*}
&& R_N \Biggl(\sum_{j=1}^M
t_j f_j \Biggr)=\frac{1}{N}\sum
_{i=1}^N \Biggl(Y_i-\sum
_{j=1}^M t_j f_j(X_i)
\Biggr)^2=\sum_{j=1}^M
\frac
{N_j}{N} (1- t_j )^2.
\end{eqnarray*}
For $\xi>0$ define $\hat t(\xi)\in\mathbb{R}^M$ by setting
\[
\hat t(\xi)_j=\cases{
1, &\quad\mbox{if there exists }$i\in\{1,\ldots,N\}$ \mbox{ s.t. } $X_i\in
\mathcal{X}_j$,\vspace*{3pt}\cr
\xi, &\quad \mbox{if there is no } $i\in\{1,\ldots,N\}$ \mbox{ s.t. } $X_i
\in\mathcal{X}_j$.}
\]
Hence, $\hat t(\xi) \in\operatorname{arg}\operatorname{min}_{t\in
\mathbb{R}^M}R_N(\sum_{j=1}^M t_j
f_j)$ and
$\hat{h}_\xi=\sum_{j=1}^M \hat t(\xi)_j f_j$ is an empirical
minimizer in $\operatorname{span}(F)$.

For every sample in $\Omega_0$, let
$j_0\in\{1,\ldots,N\}$ be the index for which $X_i \notin\mathcal{X}
_{j_0}$ for every $1 \leq i \leq N$. Therefore,
\[
R(\hat{h}_\xi)=\mathbb{E} \bigl(Y-\hat{h}_\xi(X)
\bigr)^2\geq(\xi-1)^2 P(X\in\mathcal{X}_{j_0})=
\frac{(\xi-1)^2}{k}
\]
and the claim follows by selecting $\xi$ large enough.

\begin{appendix}
\section*{Appendix}
\label{secappendix}

We begin by presenting a proof of the well-known fact that if the
$L_\infty$ and $L_2$ norms are $\sqrt{B}$-equivalent on the span of
$M$ linearly-independent functions, then $B \geq M$.

Let $F=\{f_1,\ldots, f_M\} \subset L_2$ be a dictionary whose span is
of dimension $M$, and recall
that
\begin{equation}
\label{defB-2}
\sqrt{B}=\sup_{f\in\operatorname{span}(F)
\setminus\{0\}}\frac{\Vert f\Vert_{L_\infty}}{\Vert f\Vert_{L_2}}.
\end{equation}
For every $u \in\mathbb{R}^M$ set $f_u=\sum_{j=1}^M u_j f_j$ and
define an
inner-product on $\mathbb{R}^M$ by
\[
\langle u,v \rangle_F=\mathbb{E}f_u(X)f_v(X).
\]
Let $(v_1,\ldots,v_M)$ be an orthonormal basis
of $\mathbb{R}^M$ relative to $ \langle\cdot,\cdot
\rangle_F$ and for every
$1\leq j\leq M$, set $\phi_j=f_{v_j}$. Observe that $(\phi_1,\ldots
,\phi_M)$ is an
orthonormal basis of $\operatorname{span}(F)$ in $L_2$.

For $\mu$-almost every $x\in\mathcal{X}$,
\[
\sum_{j=1}^M \phi_j^2(x)
\leq\mathop{\operatorname{ess}\operatorname{sup}}_{z\in\mathcal{X}}\sum
_{j=1}^M \phi_j(x)\phi_j(z)=
\Biggl\Vert\sum_{j=1}^M \phi_j(x)
\phi_j\Biggr\Vert _{L_\infty},
\]
and by the definition of $B$ in
(\ref{defB-2}),
\[
\Biggl\Vert\sum_{j=1}^M \phi_j(x)
\phi_j\Biggl\Vert_{L_\infty}\leq\sqrt {B}\Biggl\Vert\sum
_{j=1}^M \phi_j(x)\phi_j
\Biggr\Vert_{L_2}= \sqrt{B} \Biggl(\sum_{j=1}^M
\phi_j^2(x) \Biggr)^{1/2}.
\]
Hence,\vspace*{-3pt} for $\mu$-almost every $x\in\mathcal{X}$,
\[
\sum_{j=1}^M \phi_j^2(x)
\leq B,
\]
and by integrating this inequality with respect to $\mu$ and recalling
that $\mathbb{E}\phi_j^2(X)=1$,
it follows that $M\leq B$.

\begin{pf*}{Proof of Proposition~\protect\ref{propexpo-moments}}
Consider the model $Y= \langle X,t^*  \rangle+\zeta$ where
$t^*\in\mathbb{R}^M$, $X$ is a standard Gaussian vector in $\mathbb
{R}^M$ and $\zeta
$ is a mean-zero noise that is independent of $X$. To make the
presentation simpler, assume that $t^*=0$, and thus one only observes
the noise $Y=\zeta$. The aim here is to estimate the distance between
$\hat{t}$ and $t^*=0$ when the
noise $\zeta$ is only assumed to be in $L_2$.

Let us begin by showing that, conditionally on
$\zeta_1,\ldots,\zeta_N$, and if $\hat\sigma_N^2=\frac{1}{N}\sum_{i=1}^N \zeta_i^2$, then with probability at least $1-2\exp(-c_0 N)$,
\begin{equation}
\label{eqmainresult1} R(\hat t)-R \bigl(t^* \bigr)= \Vert\hat t\Vert_2^2
\geq\frac{c\hat\sigma_N^2 M}{N},
\end{equation}
for a suitable absolute constant $c$.

To that end, observe that the excess empirical risk for every $v\in
\mathbb{R}
^M$ is
\begin{equation}
\label{eqdecomp-q-ll} P_N\mathcal{L}_v=R_N(v)-R_N(0)=
\frac{1}{N}\sum_{i=1}^N \langle
X_i,v \rangle^2-\frac
{2}{N}\sum
_{i=1}^N \zeta_i \langle
X_i,v \rangle,
\end{equation}
and that for every sample, if $r_1 < r_2$ and
\[
\inf_{0\leq r< r_1}\inf_{\Vert v\Vert_2=r} P_N
\mathcal{L}_v> \inf_{r\geq
r_2}\inf_{\Vert v\Vert_2=r}
P_N\mathcal{L}_v,
\]
one\vadjust{\goodbreak} has $\|\hat{t}\|_2 \geq r_1$.

Using a standard $\varepsilon$-net argument together with Gaussian
concentration, one may show that if $N\geq c_0 M$, then with $\mu
^N$-probability at least
$1-2\exp(-c_1 N)$, for every $x\in\mathbb{R}^M$,
\begin{equation}
\label{eqrip} \frac{1}{2}\Vert x\Vert_2^2\leq
\frac{1}{N}\sum_{i=1}^N \langle
X_i,x \rangle^2\leq\frac{3}{2}\Vert x
\Vert_2^2.
\end{equation}
Moreover, on that event, setting
\[
I = \sup_{\{x\in\mathbb{R}^M \dvt\Vert x\Vert_2=1\}} \Biggl\vert \frac{1}{N}\sum
_{i=1}^N \zeta_i \langle
X_i,x \rangle \Biggr\vert,
\]
one has that for any $\zeta_1,\ldots,\zeta_N$
\[
c_1\hat\sigma_N\sqrt{\frac{M}{N}} \leq I \leq
c_2\hat\sigma_N\sqrt{\frac{M}{N}}
\]
for suitable absolute constants $c_1$ and $c_2$. We refer the reader to
Lemma~2.6.4 and Theorem~2.6.5 in
\cite{MR3113826} for more details on the techniques used to obtain
these observations.

Clearly, for every $r>0$,
\begin{equation}
\label{eqlinear} \inf_{\{x\in\mathbb{R}^M \dvt\Vert x\Vert_2=r\}} \frac{1}{N}\sum
_{i=1}^N \zeta_i \langle
X_i,x \rangle = -rI.
\end{equation}
Hence, by (\ref{eqdecomp-q-ll}), it follows that for $N \geq c_0M$
and conditioned on $\zeta_1,\ldots,\zeta_N$, with probability at
least $1-2\exp(-c_3N)$,
\begin{eqnarray*}
&&\inf_{0\leq r < I/6}\inf_{\Vert v\Vert_2=r}P_N
\mathcal{L}_v\geq\inf_{0\leq
r< I/6} \biggl(
\frac{r^2}{2}-rI \biggr)
\\
&&\quad > \inf_{r\geq I/3} \biggl(\frac{3r^2}{2}-rI \biggr)\geq\inf
_{r \geq I/3}\inf_{\Vert v\Vert_2=r}P_N
\mathcal{L}_v.
\end{eqnarray*}
Therefore, on that event
\[
\Vert\hat t\Vert_2\geq I/6 \geq c_4\hat
\sigma_N \sqrt{\frac{M}{N}}.
\]

Now, all that remains is to show that $P(\hat{\sigma}_N^2 \geq x)
\geq c_5/x$.
\end{pf*}

\begin{Lemma}\label{lemlower-bound-noise}
For every $N\geq2$ and $x\geq1$, there exists
a mean-zero, variance one random variable $\zeta$ for which
\[
P \bigl(\hat\sigma_N^2\geq x \bigr)\geq
\frac{c_1}{x}.
\]
\end{Lemma}

\begin{pf}
Fix $x \geq1$, let $\varepsilon$ be a symmetric, $\{-1,1\}$-valued
random variable, set $\delta=1/(xN)$ and put $\eta$ to be a $\{0,1\}
$-valued random variable with mean $\delta$ that is independent of
$\varepsilon$. Finally, let $R=1/\sqrt{\delta}$ and set $\zeta
=R\varepsilon\eta
$. Thus, $\mathbb{E}\zeta=0$ and $\Vert\zeta\Vert_{L_2}= R \delta
^{1/2}=1$.

Let $\zeta_i=R\varepsilon_i\eta_i$, $i=1,\ldots,
N$ be independent copies of $\zeta$. Recall that $NR^{-2}x=1$ and that
$\delta N \leq1$. Therefore,
\begin{eqnarray*}
&&P \bigl(\hat\sigma_N^2\geq x \bigr)=P \Biggl(
\frac{1}{N}\sum_{i=1}^N
\zeta_i^2\geq x \Biggr)= P \Biggl(\sum
_{i=1}^N \eta_i\geq1 \Biggr)
\\[-2pt]
&&\quad= P \bigl(\exists i\in\{1,\ldots,N\}, \eta_i=1 \bigr)=1-(1-
\delta)^N\geq c_1N\delta=c_1/x,
\end{eqnarray*}
as claimed.\vspace*{-9pt}
\end{pf}
\end{appendix}

\section*{Acknowledgment}\vspace*{-3pt}
Shahar Mendelson was supported by the Mathematical Sciences Institute
-- The Australian National University and by ISF Grant\vspace*{-8pt} 900/10.

%

%

%



\printhistory
\end{document}